\documentclass[12pt]{amsart}
\usepackage{geometry}                
\usepackage{graphicx}
\usepackage{amssymb}
\usepackage{amscd}
\usepackage{epstopdf}
\usepackage{enumerate}
\newtheorem{df}{Definition}[section]

\newtheorem{thm}{Theorem}[section]
\newtheorem{prop}{Proposition}[section]
\newtheorem{lm}{Lemma}[section]

\newtheorem{remark}{Remark}[section]
\newtheorem{fact}{Fact}[section]
\newtheorem{cor}{Corollary}[section]
\title{A tower of Ramanujan graphs and a reciprocity law of graph zeta functions}
\author{Kennichi Sugiyama}

\begin{document}
\maketitle

\begin{center}
Department of Mathematics, Faculty of Science,\\ 
Rikkyo University, 3-34-1 Nishi-Ikebukuro, Toshima,\\
Tokyo 171-8501, Japan \\
e-mail address : kensugiyama@rikkyo.ac.jp
\end{center}

\begin{abstract} 
Let $l$ be an odd prime. We will construct a tower of connected regular Ramanujan graph of degree $l+1$ from of modular curves. This supplies an example of a collection of graphs whose discrete Cheeger constants are bounded by $(\sqrt{l}-1)^{2}/2$ from below. We also show graph (or Ihara) zeta functions satisfy a certain reciprocity law. \\

AMS classification 2010: 05C25, 05C38, 05C50, 05C75, 11G18, 11G20, 11M38, 11M99.
\end{abstract}

\section{Introduction}

Let $p$ be a prime satisfying $p\equiv 1\,({\rm mod}\,12)$ and let us fix an odd prime $l$ different from $p$. In \cite{Sugiyama2017} we have constructed a connected regular Ramanujan graph $G^{(l)}_p(1)$ of degree $l+1$ non-bipartite.  The number of vertices $G^{(l)}_p(1)$ is $(p-1)/12$ and the Euler characteristic is
\[\chi(G^{(l)}_p(1))=\frac{(p-1)(1-l)}{24}.\]
The graph $G^{(l)}_p(1)$ is regarded as a graph of level {\em one}. In this paper we will construct a connected non-bipartite regular Ramanujan graph of degree $l+1$ of a higher level.\\

In the following let $p$ be a prime such that $p\equiv 1\,({\rm mod}\,12)$ and $l$ an odd prime different from $p$. Let ${\mathcal N}_{p,l}$ be the set of square free positive integers such that every member $N$ is prime to $lp$. Then to each $N$ of ${\mathcal N}_{p,l}$, a connected non-bipartite regular Ramanujan graph $G^{(l)}_p(N)$ of degree $l+1$ will be assigned. 
  Let $\lambda_0(G^{(l)}_p(N))\leq \lambda_1(G^{(l)}_p(N)) \leq \cdots \leq \lambda_{\nu(N)-1}(G^{(l)}_p(N))$ denote eigenvalues of the Laplacian of $G^{(l)}_p(N)$. Since $G^{(l)}_p(N)$ is connected $\lambda_0(G^{(l)}_p(N))=0$ and $\lambda_1(G^{(l)}_p(N))$ is positive. A relationship between the adjacency matrix and the Laplacian (cf. (2)) shows that
 \begin{equation}\rho^{i}(G^{(l)}_p(N)):=(l+1)-\lambda_i(G^{(l)}_p(N))\end{equation}
 is an eigenvalue of the adjacency matrix. 
%
\begin{thm} 
\begin{enumerate} 
\item For $i \geq 1$.
\[(\sqrt{l}-1)^2 \leq \lambda_i(G^{(l)}_p(N)) \leq (\sqrt{l}+1)^2,\quad \forall N \in {\mathcal N}_{p,l}.\]

\item Let $M$ and $N$ be elements of ${\mathcal N}_{p,l}$ satisfying $M|N$. Then $G^{(l)}_p(N)$ is a covering of $G^{(l)}_p(M)$ of degree $\sigma_1(N/M)$ and
\[\rho^1(G^{(l)}_p(N)) \geq \rho^1(G^{(l)}_p(M)),\quad \lambda_1(G^{(l)}_p(N)) \leq \lambda_1(G^{(l)}_p(M)).\]
Here $\sigma_1$ is the Euler function defined by
\[\sigma_1(n)=\sum_{d|n}d.\]
\end{enumerate}
\end{thm}
Our tower of Ramanujan graphs $\{G^{(l)}_p(N)\}_{N\in {\mathcal N}_{p,l}}$ has an interesting geometric property. In order to explain further we recall {\em the (discrete) Cheeger constant}. In general let $G$ be a connected $d$-regular graph of $n$ vertices. The Cheeger constant $h(G)$ of $G$ is defined by
\[h(G)={\rm min}\{\frac{|\partial S|}{|S|}\,:\, S \subset V(G), \, 0<|S|\leq \frac{n}{2}\},\]
where $V(G)$ denotes the set of vertices and
\[\partial S:=\{\{u,v\}\in GE(G)\,:\, u\in S,\,\,v\in V(G)\setminus S\}.\]
Here $GE(G)$ is the set of geometric edges (i.e. the set of unoriented edges, see \S 2) and $|\cdot|$ denotes the cardinality.
Then the smallest non-zero eigenvalue $\lambda_1(G)$ of the Laplacian satisfies (\cite{Alon-Milman} \cite{Tanner})
\[\frac{\lambda_1(G)}{2} \leq h(G) \leq \sqrt{2d\lambda_1(G)}\]
and the next corollary is an immediate consequence of {\bf Theorem 1.1}.
\begin{cor} (A gap theorem)
\[\frac{(\sqrt {l} -1)^2}{2} \leq h(G^{(l)}_{p}(N)) \leq \sqrt{2(l+1)}(\sqrt{l}+1)\]
for any $N\in {\mathcal N}_{p,l}$.
\end{cor}
In general the graph zeta function (or the Ihara zeta function) $Z(G)(t)$ is defined for a finite connected graph $G$. Although a priori $Z(G)(t)$ is a power series of $t$, the Ihara formula tells us that it is a rational function (see {\bf Fact 2.1}). We will show that the zeta functions of our graphs satisfy a reciprocity law.
%
%
%
\begin{thm} (A reciprocity law)
Let $p$ and $q$ be distinct primes satisfying $p\equiv q \equiv 1\,({\rm mod}\,12)$ and $l$ an odd prime different from $p$ and $q$. Then 
\[\frac{Z(G^{(l)}_p(q))(t)}{Z(G^{(l)}_p(1))(t)^2}=\frac{Z(G^{(l)}_q(p))(t)}{Z(G^{(l)}_q(1))(t)^2}.\]
In particular
\[Z(G^{(l)}_p(q))(t)\equiv Z(G^{(l)}_q(p))(t)\quad {\rm mod}\,{\mathbb Q}(t)^{\times 2}.\]
\end{thm}
Here is an application of {\bf Theorem 1.1} to modular forms.  As before let $p$ be a prime satisfying $p\equiv 1\,({\rm mod}\,12)$ and $N$ a square free positive integer prime to $p$. Then the spaces of cusp forms $S_2(\Gamma_0(pN))$ and one of $p$-new forms $S_2(\Gamma_0(pN))_{pN/N}$ of level $pN$ (see \S 4, especially (21)) have decompositions
\[S_2(\Gamma_0(pN))=\oplus_{\alpha}{\mathbb C}f_{\alpha},\quad S_2(\Gamma_0(pN))_{pN/N}=\oplus_{\chi}{\mathbb C}f_{\chi},\]
where $f_{\alpha}$ and $f_{\chi}$ are normalized Hecke eigenforms of character $\alpha$ and $\chi$ (cf. {\bf Theorem 4.1} and (22)). 
%
Using the result due to Alon-Boppana (\cite{Alon} \cite{Alon-Milman}) we will show the  following.
\begin{thm} Let $p$ be a prime satisfying $p\equiv 1\,({\rm mod}\,12)$ and $l$ an odd prime different from $p$. Let $\{r_i\}_{i=1}^{\infty}$ be a set of mutually distinct primes not dividing $lp$. Set $N_k=\prod_{i=1}^{k}r_i$ and then
\[\lim_{k\to \infty}{\rm Max}\{a_l(f_{\chi})\,:\, S_2(\Gamma_0(pN_k))_{pN_k/N_k}=\oplus_{\chi}{\mathbb C}f_{\chi}\}=2\sqrt{l},\]
where $a_l(f_{\chi})$ denotes the $l$-th Fourier coefficient of $f_{\chi}$. In particular
\[\lim_{k\to \infty}{\rm Max}\{a_l(f_{\alpha})\,:\, S_2(\Gamma_0(pN_k))=\oplus_{\alpha}{\mathbb C}f_{\alpha}\}=2\sqrt{l}.\]

\end{thm}

\section{Basic facts of the zeta function of a graph}
A (finite) graph $G$ consists of a finite set of vertices $V(G)$ and a finite set of oriented edges $E(G)$, which satisfy the following property: there are {\it end point maps},
\[\partial_0, \quad \partial_1 : E(G) \to V(G),\]
and {\it an orientation resersal},
\[J:E(G) \to V(G),\quad J^2=\text{identity},\]
such that $\partial_i\circ J=\partial_{1-i}\,(i=0,1)$. The quotient $E(G)/J$ is called {\it the set of geometric edges} and is denoted by $GE(G)$. We regard an element of $e\in GE(G)$ as an unoriented edge and if its end-points are $u$ and $v$ we write $e=\{u,v\}$.  For $x\in V(G)$ we set
\[E_j(x)=\{e\in E(G)\, | \, \partial_j(e)=x\}, \quad j=0,1.\]
Thus $JE_j(x)=E_{1-j}(x)$. Intuitively $E_0(x)$ (resp. $E_1(x)$) is the set of edges departing from (resp. arriving at) $x$. The {\it degree} of $x$, $d(x)$, is defined by
\[d(x)=|E_0(x)|=|E_1(x)|.\]
$E(G)$ is naturally divided into two classes, {\it loops} and {\it  passes}. An edge $e\in E(G)$ is called {\it a loop} if $\partial_0(e)=
\partial_1(e)$ and is called {\it a pass} otherwise. Let $2l(x)$ and $p(x)$ be the number of loops and passes starting from $x$, respectively (both $l(x)$ and $p(x)$ are positive integers). Note that, because of the involution $J$, if we replace "departing" by "arriving" these number does not change. By definition, it is clear that
\[d(x)=2l(x)+p(x).\]
We set $q(x):=d(x)-1$. 
Let $C_{0}(G)$  be the free ${\mathbb Z}$-module generated by $V(G)$ with vertices as the natural basis. We define endomorphisms $Q$ and $A$ of $C_0(G)$ by
\[\quad Q(x)=q(x)x, \quad x\in V(G),\]
and 
\[A(x)=\sum_{e\in E(G), \partial_0(e)=x}\partial_1(e),\quad x\in V(G),\]
respectively. Note that because of the involution $J$, 
\[A(x)=\sum_{e\in E(G), \partial_1(e)=x}\partial_0(e).\]
The operator $A$ will be called {\it the adjacency operator}.  We sometimes identify it with the representing matrix with respect to the basis $\{x\}_{x\in V(G)}$. Thus the $yx$-entry $A_{yx}$ of $A$ is the number of edges departing from $x$ and arriving at $y$. The orientation reversing involution $J$ implies
\[A_{xy}=A_{yx}. \]
Note that $A_{xx}=2l(x)$ and $p(x)=\sum_{y\neq x}A_{yx}$. If $d(x)=k$ for all $x\in V(G)$, $G$ is called $k$-{\it regular}.\\

Connecting distinct vertices $x$ and $y$ by geometric $A_{xy}$-edges and drawing $\frac{1}{2}A_{xx}$-loops at $x$,  the adjacency matrix $A$ determines an unoriented $1$-dimensional simplicial complex. We call it {\it the geometric realization} of $G$, and denote it by $G$ again. We say that $G$ is connected if the geometric realization is. The Euler characteristic $\chi(G)$ is equal to $|V(G)|-|GE(G)|$, hence if $G$ is connected, the fundamental group is a free group of rank $1-|V(G)|+|GE(G)|$. 
For a later purpose, we summarize the relationship between a graph and its adjacency matrix.
\begin{prop} Let $A=(a_{ij})_{1\leq i,j \leq m}$ be an $m\times m$-matrix satisfying the following conditions.
\begin{enumerate}[(1)]
\item The entries $\{a_{ij}\}_{ij}$ are non-negative integers and satisfy
\[a_{ij}=a_{ji}, \quad \forall i \,\text{and}\, j.\]
\item $a_{ii}$ is even for every $i$.
\end{enumerate}
Then there is a unique graph $G$ whose adjacency matrix is $A$. Moreover, $G$ is $k$-regular if and only if one of the following equivalent condition satisfied :
\begin{enumerate}[(a)]
\item
\[\sum_{i=1}^{m}a_{ij}=k,\quad \forall j\]
\item 
\[\sum_{j=1}^{m}a_{ij}=k,\quad \forall i.\]
\end{enumerate}
\end{prop}
 In the following, a graph $G$ is always assumed to be {\it connected}. {\it A path of length} $m$ is a sequence $c=(e_1,\cdots,e_m)$ of edges such that $\partial_0(e_i)=\partial_1(e_{i-1})$ for all $1< i \leq m$ and the path is {\it reduced} if $e_i\neq J(e_{i-1})$ for all      $1< i \leq m$. The path is {\it closed} if $\partial_0(e_1)=\partial_1(e_{m})$, and the closed path has {\it no tail} if $e_m\neq J(e_1)$. A closed path of length one is nothing but a loop. Two closed paths are {\it equivalent} if one is obtained from the other by a cyclic shift of the edges. Let ${\frak C}(G)$ be the set of equivalence classes of reduced and tail-less closed paths of $G$. Since the length depends only on the equivalence class, the length function descends to the map;
\[l : {\frak C}(G) \to {\mathbb N},\quad l([c])=l(c),\]
where $[c]$ is the class determined by $c$. We define a reduced and tail-less closed path $C$ to be primitive if it is not obtained by going $r\,(\geq 2)$ times some another closed path. Let ${\frak P}(G)$ be the subset of ${\frak C}(G)$ consisting of the classes of primitive closed paths (which are  reduced and tail-less by definition). The graph zeta function ( or {\em Ihara zeta function}) of $G$ is defined to be
\[Z(G)(t)=\prod_{[c]\in {\frak P}(G)}\frac{1}{1-t^{l([c])}}.\]
Although this is an infinite product, it is a rational function.
\begin{fact}(\cite{Bass},\cite{Hoffman},\cite{Ihara1966},\cite{Stark-Terras})
\[Z(G)(t)=\frac{(1-t^2)^{\chi(G)}}{{\rm det}[1-At+Qt^2]}.\]
\end{fact}
\begin{fact}(\cite{Sugiyama2017})
Let $G$ be a $k$-regular graph with $m$ vertices. Then the Euler characteristic $\chi(G)$ is
\[\chi(G)=\frac{m(2-k)}{2}.\]
\end{fact}
\begin{remark}
Note that the Euler characteristic does not depend on the number of loops.
\end{remark}
Let $E_{or}(G)\subset E(G)$ be a section of the natural projection $E(G) \to GE(G)$. In other word we choose an orientation on geometric edges and make the geometric realization into an oriented one dimensional simplicial complex. Let $C_1(G)$ be the free ${\mathbb Z}$-module generated by $E_{or}(G)$. Then the boundary map 
\[\partial : C_1(G) \to C_0(G)\]
is naturally defined. Let $\partial^{t}$ be the dual of $\partial$ and the {\it Laplacian} $\Delta$ of $G$ is defined to be $\Delta=\partial \partial^{t}$. It is known (and easy to check) that (\cite{Terras}, \cite{Hoffman}),
\begin{equation}\Delta=1-A+Q.\end{equation}
Now let $G$ be a connected $k$-regular graph. Since $0$ is an eigenvalue of $\Delta$ with multiplicity one, (2) shows that $k$ is an eigenvalue of $A$ with multiplicity one. Because of semi-positivity of $\Delta$ we find that 
\[|\lambda|\leq k\quad \text{for any eigenvalue $\lambda$ of $A$}\]
and that $-k$ is an eigenvalue of $A$ if and only if $G$ is bipartite (\cite{Terras}, {\bf Chapter 3}). Here $G$ is called {\em bipartite} if the set of vertices $V(G)$ can be divided into disjoint subset $V_0$ and $V_1$ such that every edge connects points in $V_0$ and $V_1$, namely there is no edge whose end points are simultaneously contained in $V_{i}$ ($i=0,1$). 

\begin{df} Let $G$ be a $k$-regular graph. We say that it is {\rm Ramanujan}, if all eigenvalues $\lambda$ of $A$ with $|\lambda|\neq k$ satisfy
\[|\lambda|\leq 2\sqrt{k-1}.\]
\end{df}
See \cite{Li1992}, \cite{Murty} and \cite{Valette} for detailed expositions of Ramanujan graphs. \\

A map $f$ from a graph $G^{\prime}$ to $G$ is defined to be a pair $f=(f_V,f_E)$ of maps 
\[f_V : V(G^{\prime}) \to V(G),\quad f_E : E(G^{\prime}) \to E(G)\]
 satisfying
\[\partial_{i}f_E=f_V\partial_{i},\quad i=0,1.\] 
Suppose that $G$ and $G^{\prime}$ are connected. If there is a positive integer $d$ such that $|f_V^{-1}(v)|=|f_E^{-1}(e)|=d$ for any $v\in V(G)$ and $e\in E(G)$, $f$ is mentioned as {\em a covering map of degree $d$}.
\section{A construction of a Ramanujan graph}
Let $p$ be a prime, and $B$ the quaternion algebra over ${\mathbb Q}$ ramified at two places $p$ and $\infty$. Let $R$ be a fixed maximal order in $B$ and $\{I_1,\cdots,I_n\}$ be the set of left $R$-ideals representing the distinct ideal classes. We choose $I_1=R$ and say $n$ {\it the class number} of $B$. For $1\leq i \leq n$, $R_i$ denotes the right order of $I_i$, and let $w_i$ the order of $R_i^{\times}/\{\pm 1\}$. 
The product 
\begin{equation}W=\prod_{i=1}^{n}w_{i}\end{equation}
is independent of the choice of $R$ and is equal to the exact denominator of $\frac{p-1}{12}$ (\cite{Gross} p.117) and Eichler's mass formula states that
\[\sum_{i=1}^{n}\frac{1}{w_i}=\frac{p-1}{12}.\]
Let ${\mathbb F}$ be an algebraic closure of ${\mathbb F}_p$. There are $n$ distinct isomorphism classes $\{E_1,\cdots,E_n\}$ of supersingular elliptic curves over ${\mathbb F}$ such that ${\rm End}(E_i)\simeq R_i$.\\
Now we assume that $p-1$ is divisible $12$. Then $\frac{p-1}{12}$ is an integer and $W=\prod_{i=1}^{n}w_{i}=1$, namely $w_i=1$ for all $i$. Hence by Eichler's mass formula 
\begin{equation}
n=\frac{p-1}{12}.
\end{equation}

We fix an odd prime $l$ different from $p$ and let ${\mathcal N}_{p,l}$ denote the set of square free positive integers prime to $lp$. For $N\in {\mathcal N}_{p,l}$, {\em an enhanced supersingular elliptic curve of level $N$} is defined to be a pair ${\mathbf E}=(E, C_N)$ of a supersingular elliptic curve $E$ and its cyclic subgroup $C_N$ of order $N$. A homomorphism $\phi$ from ${\mathbf E}=(E, C_N)$ to ${\mathbf E}^{\prime}=(E^{\prime}, C^{\prime}_N)$ is defined by a homomorphism $\phi : E \to E^{\prime}$ satisfying
\[\phi(C_N)=C^{\prime}_N.\]
Let $\Sigma_{N}$ be the set of isomorphism classes of enhanced supersingular elliptic curve of level $N$ defined over ${\mathbb F}$. Then the cardinality $\nu(N)$ of $\Sigma_{N}$ is
\begin{equation}\nu(N)=\frac{(p-1)\sigma_{1}(N)}{12}, \quad \sigma_1(N)=\sum_{d\mid N}d.\end{equation}
Here $\sigma_1(N)$ counts the number of cyclic subgroups of $E$ of order $N$. Let ${\rm Hom}({\mathbf E}_i,\,{\mathbf E}_j)(l)$ denote the set of homomorphisms from ${\mathbf E}_i$ to ${\mathbf E}_j$ of degree $l$. We define {\em the Brandt matrix} $B_p^{(l)}(N)$ is defined to be a $\nu(N)\times \nu(N)$-matrix whose $(i,j)$-entry is 
\begin{equation}b_{ij}=\frac{1}{2}|{\rm Hom}({\mathbf E}_j,\,{\mathbf E}_i)(l)|.\end{equation}

\begin{prop} Let $N\in {\mathcal N}_{p,l}$. Then the Brandt matrix $B_p^{(l)}(N)=(b_{ij})_{1\leq i,j \leq \nu(N)}$ satisfies the following.
\begin{enumerate}[(1)]
\item Every entry is a non-negative integer and $B_p^{(l)}(N)$ is symmetric;
\[b_{ij}=b_{ji}.\]
\item
The diagonal entires $\{b_{ii}\}_{i}$ are even for all $i$.
\item For any $i=1,\cdots, \nu(N)$,
\[\sum_{j=1}^{n}b_{ij}=l+1.\]
\end{enumerate}
\end{prop}
{\bf Proof.} By definition a homomorphism from ${\mathbf E}_i=(E_i,C_N)$ to ${\mathbf E}_j=(E_j,D_N)$ is a homomorphism $\phi : E_i \to E_j$ of degree $l$ satisfying
\[\phi(C_N)=D_N.\]
%
Being $\check{\phi}$ the dual of $\phi$, $\check{\phi}\phi=l$ and $\check{\phi}(D_N)=\check{\phi}(\phi(C_N))=C_N$. Hence taking the dual homomorphisms yields bijective correspondence
\[I : {\rm Hom}({\mathbf E}_i,\,{\mathbf E}_j)(l) \to {\rm Hom}({\mathbf E}_j,\,{\mathbf E}_i)(l),\quad I(\phi)=\check{\phi},\]
which implies (1). In order to show the claim (2), it is sufficient to show that the action of $I$ on ${\rm End}({\mathbf E}_i)(l)/{\pm 1}$ has no fixed point. Let $\phi$ be an element of ${\rm End}({\mathbf E}_i)(l)/{\pm 1}$. 
Then ${\rm Ker}\phi\simeq {\rm Ker}\check{\phi}\simeq {\mathbb F}_l$ and there is a skew-symmetric nondegenerate pairing derived from the Weil paring (\cite{Silverman} \S III {\bf Remark 8.4})
\[{\rm Ker}\phi \times {\rm Ker}\check{\phi} \to \mu_{l}.\]
Suppose that there were $\phi\in {\rm End}({\mathbf E}_i)(l)/{\pm 1}$ fixed by $I$. Then $\check{\phi}=\pm \phi$ and ${\rm Ker}\phi= {\rm Ker}\check{\phi}$, which contradicts to non-degeneracy of the pairing. The claim (3) follows from the following observation : Let $E_j$ be the underlying supersingular elliptic curve of ${\mathbf E}_j$. Then by definition $\sum_{i=1}^{n}b_{ij}$ is equal to the number of cyclic subgroups of $E_i$ of order $l$, which is $l+1$.
\begin{flushright}
$\Box$
\end{flushright}
By {\bf Proposition 2.1} there is a regular graph $G_p^{(l)}(N)$ of degree $l+1$ whose adjacency matrix is $B_p^{(l)}(N)$. In {\bf Theorem 5.1} we will show that it is a connected non-bipartite Ramanujan graph. 
\begin{thm} Let $M$ and $N$ be elements of ${\mathcal N}_{p,l}$ such that $M$ is a divisor of $N$. Then there is a covering map
\[\pi_{N/M}\,:\, G^{(l)}_{p}(N) \to G^{(l)}_{p}(M)\]
of degree $\sigma_1(N/M)$
\end{thm}
{\bf Proof.}
Since $N$ is square free $M$ and $N/M$ are coprime. Thus a cyclic subgroup $C_N$ is written by
\[C_M=C_M\oplus C_{N/M}\]
and we define
\[(\pi_{N/M})_V\,:\, V(G^{(l)}_{p}(N)) \to V(G^{(l)}_{p}(M)),\quad (\pi_{N/M})_V(E, C_M\oplus C_{N/M})=(E, C_M).\]
Since the number of cyclic subgroups of $E$ of order $N/M$ is $\sigma_1(N/M)$, $|\pi_{N/M}^{-1}(v)|=\sigma_1(N/M)$ for any $v\in V(G^{(l)}_{p}(M))$. By definition an edge of $G^{(l)}_{p}(N)$ from ${\mathbf E}=(E, C_M\oplus C_{N/M})$ to ${\mathbf E}^{\prime}=(E^{\prime}, C^{\prime}_M\oplus C^{\prime}_{N/M})$ is a homomorphism $f$ from $E$ to $E^{\prime}$
satisfying
\[f(C_M)=C^{\prime}_M, \quad f(C_{N/M})=C^{\prime}_{N/M}.\]
 Forget the homomorphism of cyclic subgroups of order $N/M$ and we have 
\[{\rm Hom}({\mathbf E}, {\mathbf E}^{\prime})(l)/\{\pm 1\} \to {\rm Hom}(\pi_{N/M}({\mathbf E}), \pi_{N/M}({\mathbf E}^{\prime}))(l)/\{\pm 1\},\]
which defines a map of the set of edges
\[(\pi_{N/M})_E\,:\, E(G^{(l)}_{p}(N)) \to E(G^{(l)}_{p}(M))\]
satisfying
\[\partial_i\circ (\pi_{N/M})_E=(\pi_{N/M})_V\circ \partial_i,\quad i=0,1.\]
One finds that this map has degree $\sigma_1(N/M)$. In fact let $g$ be an element of ${\rm Hom}(\pi_{N/M}({\mathbf E}), \pi_{N/M}({\mathbf E}^{\prime}))(l)$. Thus $g$ is a homomorphism from $E$ to $E^{\prime}$ of degree $l$ satisfying 
\[g(C_M)=C^{\prime}_M.\]
Let $C_{N/M}$ be a cyclic subgroup of $E$ of order $N/M$ and we set $C_{N/M}^{\prime}=g(C_{N/M})$. Then we have a homomorphism of enhanced supersingular elliptic curve of level $N$
\[g : (E, C_M\oplus C_{N/M}) \to (E^{\prime}, C^{\prime}_M\oplus C^{\prime}_{N/M})\]
which defines an edge of $G^{(l)}_{p}(N)$. The number of cyclic subgroups of order $N/M$ (i.e. choices of $C_{N/M}$) is $\sigma_1(N/M)$ and the claim has been proved.

\begin{flushright}
$\Box$
\end{flushright}
\section{A spectral decomposition of the character group}
For a positive integer $N$, let $S_2(\Gamma_0(N))$ denote the space of cusp forms of weight $2$ for the Hecke congruence subgroup
\[\Gamma_{0}(N):=\{\left(
  \begin{array}{cc}
    a & b \\ 
    c & d \\ 
  \end{array}\right) \in {\rm SL}_2({\mathbb Z}) : \quad c \equiv 0\, ({\rm mod}\,N)\}.
\]
Let $Y_{0}(N)$ be the modular curve which parametrizes isomorphism classes of a pair ${\bf E}=(E,C_N)$ of an elliptic curve $E$ and its cyclic subgroup $C_N$ of order $N$. It is a smooth curve defined over ${\mathbb Q}$ and the set of ${\mathbb C}$-valued points is the quotient of the upper half plane by 
$\Gamma_0(N)$. Let $X_0(N)$ be the compactification of $Y_0(N)$. It is a smooth projective curve defined over ${\mathbb Q}$ and has the canonical model over ${\mathbb Z}$ which has been studied by \cite{Deligne-Rapoport} and \cite{Katz-Mazur} in detail. The space of cusp forms $S_2(\Gamma_0(N))$ is naturally identified with the space of holomorphic $1$-forms $H^0(X_0(N),\Omega)$ and in particular  with the cotangent space ${\rm Cot}_0(J_0(N))$ at the origin of the Jacobian variety $J_0(N)$ of  $X_0(N)$.\\

For a prime $p$ with $(p,N)=1$, $X_0(N)$ furnishes the $p$-th Hecke operator defined by
 \begin{equation}T_p(E,C_N):=\sum_{C}(E/C, (C_N + C)/C),\end{equation}
where $C$ runs through all cyclic subgroup schemes of $E$ of order $p$. If $p$ is a prime divisor of $N$, an operator $U_p$ is defined by
 \begin{equation}U_p(E,C_N):=\sum_{C\neq D}(E/C, (C_N + C)/C)\end{equation}
 where $D$ is the cyclic subgroup of $C_N$ of order $p$. By the functoriality, Hecke operators act on $J_0(N)$ and ${\rm Cot}_0(J_0(N))=S_2(\Gamma_0(N))$ and the resulting action coincides with the usual one on $S_2(\Gamma_0(N))$ (see  \cite{Shimura1971}). We define the Hecke algebra as ${\mathbb T}_0(N):={\mathbb Z}[\{T_p\}_{(p,N)=1}, \{U_p\}_{p|N}\}]$, which is a commutative subring of ${\rm End}J_0(N)$. The effects of $T_p$ and $U_p$ on $f=\sum_{n=1}^{\infty}a_nq^n\in S_2(\Gamma_0(N))$ are
\begin{equation}f|U_p=\sum_{n=1}^{\infty}a_{pn}q^n\end{equation}
and
\begin{equation}f|T_p=\sum_{n=1}^{\infty}(a_{pn}+pa_{n/p})q^n.\end{equation}
Here $a_{n/p}=0$ if $n/p$ is not an integer. 
\begin{df} For a positive integer $M$, we define a subalgebra ${\mathbb T}_0(N)^{(M)}$ of ${\mathbb T}_0(N)$ to be the omitting of Hecke operators from ${\mathbb T}_0(N)$ whose indices are prime divisors of $M$, that is
\[{\mathbb T}_0(N)^{(M)}={\mathbb Z}[[\{T_p\}_{(p,NM)=1}, \{U_p\}_{p|N, (p,M)=1}\}.\]
We call an algebraic homomorphism from ${\mathbb T}_0(N)^{(M)}$ to ${\mathbb C}$ as {\em a character}. If the image is contained in ${\mathbb R}$ it is referred as {\em real}.
\end{df}

Let $M$ be a positive integer and $f$ an element of $S_2(\Gamma_0(M))$. For a positive integer $d$ we set
\[f^{(d)}(z)=f(dz)\in S_2(\Gamma_0(dM)).\]

\begin{df} Let $N$ be a square free positive integer and $M$ a divisor of $N$. For a divisor $d$ of $N/M$ we define
\[S_2(\Gamma_0(M))^{(d)}=\{f^{(d)}(z)\,|\, f\in S_2(\Gamma_0(M))\}\subset S_2(\Gamma_0(N)).\]
The space of old forms of level $N$ is defined to be
\[S_2(\Gamma_0(N))_{old}=\sum_{M|N, M\neq N}\sum_{d|(N/M)}S_2(\Gamma_0(M))^{(d)}\subset S_2(\Gamma_0(N))\]
and the orthogonal complement of $S_2(\Gamma_0(N))_{old}$ for the Petersson product is called by the space of {\em new forms} and denoted by $S_2(\Gamma_0(N))_{new}$. 
\end{df}
Let $N$ be a square free positive integer and $q$ a prime not dividing $N$. Since the action of $T_q$ on $S_2(\Gamma_0(N))$ is self-adjoint for the Petersson product and since $S_2(\Gamma_0(N))_{old}$ is stable by $T_q$, $S_2(\Gamma_0(N))_{new}$ is stable by ${\mathbb T}_0(N)^{(N)}$. This implies that $S_2(\Gamma_0(N))_{new}$ admits a spectral decomposition by ${\mathbb T}_0(N)^{(N)}$. We will show that $S_2(\Gamma_0(N))$ has an irreducible decomposition of multiplicity one by the action of the full Hecke algebra ${\mathbb T}_0(N)$ (cf. {\bf Theorem 4.1}). In proving the theorem, a key fact is the following,  which is mentioned as {\em multiplicity one} (\cite{Atkin-Lehner} \cite{Li}). 
\begin{fact} Let $N$ be a positive integer (which may not be square free) and $f=\sum_{n=1}^{\infty}a_nq^n$ an element of $S_2(\Gamma_0(N))$. Suppose that $a_n=0$ for all $n$ with $(n,t)=1$, where $t$ is a fixed positive integer. Then $f\in S_2(\Gamma_0(N))_{old}$.
\end{fact}
This fact shows that the above  spectral decomposition of $S_2(\Gamma_0(N))_{new}$ by ${\mathbb T}_0(N)^{(N)}$ has multiplicity one. One finds that this yields an irreducible decomposition of $S_2(\Gamma_0(N))_{new}$ for the full Hecke algebra. In fact let $f \in S_2(\Gamma_0(N))_{new}$ be the normalized eigenform of ${\mathbb T}_0(N)^{(N)}$ and $p$ a prime not dividing $N$. Since $T_p$ is selfadjoint for the Petersson product its eigenvalue is real number. Moreover $f$ is automatically a Hecke eigenform of the full Hecke algebra by the following reason. Let $\alpha$ be the character of ${\mathbb T}_0(N)^{(N)}$ associated to $f$ and $q$ be a prime divisor of $N$. Since ${\mathbb T}_0(N)$ is commutative $f|U_q$ is also a Hecke eigenform of ${\mathbb T}_0(N)^{(N)}$ whose character is $\alpha$. By the multiplicity one, $f|U_q$ should be a multiple of $f$;
\[f|U_q=\alpha_q f.\]
Defining $\alpha(U_q)=\alpha_q$,  we have a character $\alpha$ of ${\mathbb T}_0(N)$ and $f$ is the normalized Hecke eigenform of character $\alpha$. Moreover since $N$ is square free $\alpha_q=\pm1$ for $q\mid N$ (\cite{Hida} {\bf Lemma 3.2}) and $\alpha$ is real character. Thus we have an irreducible decomposition as a ${\mathbb T}_0(N)$-module
\[S_2(\Gamma_0(N))_{new}=\oplus_{\alpha}S_2(\Gamma_0(N))_{new}(\alpha)\]
by real characters and every irreducible component has dimension one. Here $S_2(\Gamma_0(N))_{new}(\alpha)$ denotes the isotypic component of $\alpha$ 
\[S_2(\Gamma_0(N))_{new}(\alpha)=\{f\in S_2(\Gamma_0(N))_{new}\,|\, f|T=\alpha(T)f,\quad \forall T\in {\mathbb T}_0(N)\},\]
which is spanned by the normalized Hecke eigenform. By the definition of the space of new forms we have
\begin{equation}S_2(\Gamma_0(N))=\oplus_{M|N}(\oplus_{d|(N/M)}S_2(\Gamma_0(M))_{new}^{(d)}).\end{equation}
Fix a divisor $M$ of $N$ and let us consider the subspace 
\[{\mathbb S}_M=\oplus_{d|(N/M)}S_2(\Gamma_0(M))_{new}^{(d)}.\]
Being $N/M=l_1\cdots l_m$ a prime decomposition, there is an isomorphism as vector spaces 
\begin{equation}{\mathbb S}_M\simeq S_2(\Gamma_0(M))_{new}^{\oplus 2^{m}}.\end{equation}
We will explicitly describe this isomorphism. 
\begin{prop} 
Let $N$ be a square free positive integer and $M$ a divisor of $N$. Let $f\in S_2(\Gamma_0(M))_{new}$ be a normalized Hecke eigenform. Then for $\epsilon=(\epsilon_{l_1},\cdots, \epsilon_{l_m})$ ($\epsilon_{l_i}=\pm$)
there is a normalized Hecke eigenform $f_{\epsilon}$ of level $N$ satisfying the following conditions.
\begin{enumerate}
\item If $q$ a prime not dividing $N/M$
\[a_q(f_{\epsilon})=a_q(f).\]
\item 
\[a_{l_i}(f_{\epsilon})=\alpha_{l_i}^{\epsilon_{l_i}}\]
where 
\[\alpha_{l_i}^{+}=\frac{a_{l_i}(f)+\sqrt{\Delta_i}}{2},\quad \alpha_{l_i}^{-}=\frac{a_{l_i}(f)-\sqrt{\Delta_i}}{2},\quad \Delta_i=a_{l_1}(f)^2-4{l_i} (<0).\]
Moreover the $2^m$ complex numbers $\{\alpha_{l_1}^{(\pm)},\cdots, \alpha_{l_m}^{(\pm)}\}$ are mutually different.
\end{enumerate}

\end{prop}
{\bf Proof.}
In general let $p$ be a prime and $F$ a square free positive integer prime to $p$. We have two degeneracy maps $\alpha_p,\,\beta_p : X_0(pF) \to X_0(F)$ defined by
\[\alpha_p(E,C_p\oplus C_F)=(E,C_F),\quad \beta_p(E,C_p\oplus C_F)=(E/C_p,(C_p\oplus C_F)/C_p),\]
which induces linear maps
\begin{equation}\alpha_p^{\ast}, \beta_p^{\ast} : S_2(\Gamma_0(F)) \to S_2(\Gamma_0(pF))\end{equation}
whose effects on $f=\sum_{n=1}^{\infty}a_nq^n\in S_2(\Gamma_0(F))$ are
\begin{equation}\alpha_p^{\ast}(f)=f=\sum_{n=1}^{\infty}a_nq^n,\quad \beta_p^{\ast}(f)=f^{(p)}=\sum_{n=1}^{\infty}a_{n}q^{pn}. \end{equation}
Let $T$ be $T_r\, (r\nmid pF)$ or $U_l\, (l\mid F)$. Then $T$ commutes with $\alpha_p$ and $\beta_p$ and 
\begin{equation}
\begin{CD}
 S_2(\Gamma_0(F))\oplus S_2(\Gamma_0(F))@>\text{$\alpha_p^{\ast}+ \beta_p^{\ast}$}>>  S_2(\Gamma_0(pF))\\
 @V\text{$(T,\,T)$}VV  @V\text{$T$}VV\\
   S_2(\Gamma_0(F))\oplus S_2(\Gamma_0(F))  @>\text{$\alpha_p^{\ast}+ \beta_p^{\ast}$}>>   S_2(\Gamma_0(pF)).
 \end{CD}
\end{equation}
Using (14) and (15) we will inductively construct $f_{\epsilon}$ by the number of prime divisors $m$. We set $M_m=Ml_1\cdots l_m\, (m\geq 1)$ and $M_0=M$. Suppose that we have constructed a desired normalized Hecke eigenform $f_{\epsilon}\in S_2(\Gamma_0(M_{m-1}))$ of character $\chi_{\epsilon}$. For a prime $r$ different from $l_m$, we let $T$ be $T_r$ or $U_r$ according to $r\nmid M_{m}$ or $r\mid M_{m-1}$, respectively. Then (15) implies
\begin{equation}
\begin{CD}
 S_2(\Gamma_0(M_{m-1}))\oplus S_2(\Gamma_0(M_{m-1}))@>\text{$\alpha_{l_{m}}^{\ast}+ \beta_{l_{m}}^{\ast}$}>>  S_2(\Gamma_0(M_{m}))\\
 @V\text{$(T,\,T)$}VV  @V\text{$T$}VV\\
   S_2(\Gamma_0(M_{m-1}))\oplus S_2(\Gamma_0(M_{m-1}))  @>\text{$\alpha_{l_{m}}^{\ast}+ \beta_{l_{m}}^{\ast}$}>>   S_2(\Gamma_0(M_{m})).
 \end{CD}
\end{equation}
Hence 
\[\alpha_{l_{m}}^{\ast}(f_{\epsilon})|T=\alpha_{l_{m}}^{\ast}(f_{\epsilon}|T)=\chi_{\epsilon}(T)\alpha_{l_{m}}^{\ast}(f_{\epsilon})\]
and
\[\beta_{l_{m}}^{\ast}(f_{\epsilon})|T=\beta_{l_{m}}^{\ast}(f_{\epsilon}|T)=\chi_{\epsilon}(T)\beta_{l_{m}}^{\ast}(f_{\epsilon}).\]
Define a character
\[\chi_{\epsilon}^{(l_m)} : {\mathbb T}_0(M_{m})^{(l_m)} \to {\mathbb C}\]
by
\[\chi_{\epsilon}^{(l_m)}(T)=\chi_{\epsilon}(T),\]
and $\alpha_{l_{m}}^{\ast}(f_{\epsilon})$ and $\beta_{l_{m}}^{\ast}(f_{\epsilon})$ are ${\mathbb T}_0(M_{m})^{(l_m)}$-eigenforms of the same character $\chi_{\epsilon}^{(l_m)}$. Let us investigate the action of $U_{l_m}$. By (9), (10) and (14)
\[\left(  \begin{array}{c}
    \alpha_{l_m}^{\ast}(f_{\epsilon})|U_{l_m} \\ 
    \beta_{l_m}^{\ast}(f_{\epsilon})|U_{l_m}\\ 
  \end{array}
\right)=\left(  \begin{array}{cc}
    a_{l_m}(f_{\epsilon}) & -l_m \\ 
    1 & 0 \\ 
  \end{array}
\right)\left(  \begin{array}{c}
    \alpha_{l_m}^{\ast}(f_{\epsilon}) \\ 
    \beta_{l_m}^{\ast}(f_{\epsilon})\\ 
  \end{array}
\right).\]
Use the assumption (1) and the characteristic polynomial of $U_{l_m}$ is
\[\Phi(t)=t^2-a_{l_m}(f_{\epsilon})t+l_m=t^2-a_{l_m}(f)t+l_m.\]
Since $f$ is a normalized ${\mathbb T}_0(M)$-eigenform which is new, the discriminant $\Delta_m=a_{l_m}(f)^2-4{l_m}$ is negative (\cite{C-E}). Therefore the eigenvalue of $U_{l_m}$ are mutually distinct and contained in ${\mathbb C}\setminus {\mathbb R}$. Set
\begin{equation}\alpha_{l_m}^{+}=\frac{a_{l_m}(f)+\sqrt{\Delta_m}}{2},\quad \alpha_{l_m}^{-}=\frac{a_{l_m}(f)-\sqrt{\Delta_m}}{2}\end{equation}
and let $f_{\epsilon}^{+}$ and $f_{\epsilon}^{-}$ be the corresponding normalized cusp form of level $M_m$ satisfying
\[f_{\epsilon}^{+}\mid U_{l_m}=\alpha_{l_m}^{+}f_{\epsilon}^{+},\quad f_{\epsilon}^{-}\mid U_{l_m}=\alpha_{l_m}^{-}f_{\epsilon}^{-}.\]
Extend $\chi_{\epsilon}^{(l_m)}$ to a character $\chi_{\epsilon}^{+}$ and $\chi_{\epsilon}^{-}$ of ${\mathbb T}_0(M_m)={\mathbb T}_0(M_m)^{(l_m)}[U_{l_m}]$ by
\[\chi_{\epsilon}^{+}(U_{l_m})=\alpha_{l_m}^{+},\quad \chi_{\epsilon}^{-}(U_{l_m})=\alpha_{l_m}^{-}.\]
Then $f_{\epsilon}^{+}$ and $f_{\epsilon}^{-}$ are ${\mathbb T}_0(M_m)$-eigenforms whose characters are $\chi_{\epsilon}^{+}$ and $\chi_{\epsilon}^{-}$, respectively. Observe that $\alpha_{l_m}^{+}$ and $\alpha_{l_m}^{-}$ are different from each of $\{\alpha_{l_i}^{+}, \alpha_{l_i}^{-}\}_{1\leq i \leq m-1}$, where
\[\alpha_{l_i}^{+}=\frac{a_{l_i}(f)+\sqrt{\Delta_i}}{2},\quad \alpha_{l_i}^{-}=\frac{a_{l_i}(f)-\sqrt{\Delta_i}}{2},\quad \Delta_i=a_{l_i}(f)^2-4{l_i}.\]
In fact if $\alpha_{l_m}^{+}=\alpha_{l_i}^{+}$ ($1\leq i \leq m-1$), comparing their real and imaginary part we conclude
\[a_{l_m}(f)=a_{l_i}(f),\quad \Delta_m=\Delta_i\]
which implies $l_m=l_i$.
%
Thus we have constructed normalized $2^m$ Hecke eigenforms of level $M_m$ from $f$ whose characters are mutually different.

\begin{flushright}
$\Box$
\end{flushright}

{\bf Proposition 4.1} yields a spectral decomposition of multiplicity one
\begin{equation}{\mathbb S}_M=\oplus_{\beta}{\mathbb C}f_{\beta}\end{equation}
where $f_{\beta}$ is the normalized Hecke eigenform of character $\beta$. Let $M^{\prime}$ be a divisor of $N$ different from $M$ and we consider the decomposition (18) for $M^{\prime}$,
\begin{equation}{\mathbb S}_{M^{\prime}}=\oplus_{\beta}{\mathbb C}f_{\beta^{\prime}}.\end{equation}
The following lemma shows that every character $\beta$ in (18) is different from each of $\beta^{\prime}$ in (19).

\begin{lm} Let $f\in S_2(\Gamma_0(N_f))_{new}$ (resp. $g\in S_2(\Gamma_0(N_g))_{new}$) be a normalized Hecke eigenform. If there is a positive integer $t$ such that
\[a_l(f)=a_l(g)\]
for any prime $l$ with $l\nmid t$, then $f=g$.
\end{lm}
{\bf Proof.} Let $K_f$ (resp. $K_g$) be the number field generated by Fourier coefficients of $f$ and (resp. $g$) over ${\mathbb Q}$ and let $K$ be the minimal extension of ${\mathbb Q}$ that contains $K_f$ and $K_g$.  We fix a prime $l$ satisfying $l\nmid N_fN_g$ and that completely splits in $K$. Corresponding to $f$ and $g$,  there are absolutely irreducible representations
\[\rho_{f,l} : {\rm Gal}(\overline{\mathbb Q}/{\mathbb Q}) \to {\rm GL}_2({\mathbb Q}_l), \quad \rho_{g,l} : {\rm Gal}(\overline{\mathbb Q}/{\mathbb Q}) \to {\rm GL}_2({\mathbb Q}_l)\]
of the conductor $N_f$ and $N_g$ respectively which satisfy
\[{\rm det}(t-\rho_{f,l}(Frob_q))=t^2-a_q(f)t+q,\quad (q,lN_f)=1\]
and
\[{\rm det}(t-\rho_{g,l}(Frob_q))=t^2-a_q(g)t+q,\quad (q,lN_g)=1.\]
(\cite{DDT} {\bf Theorem 3.1}). Here $Frob_q$ is the Frobenius at a prime $q$. Let $S$ be a finite set of primes. Since a semi-simple representation $\rho_l : {\rm Gal}(\overline{\mathbb Q}/{\mathbb Q}) \to {\rm GL}_2({\mathbb Q}_l)$ is determined by values ${\rm Tr}\rho_l (Frob_q)$ on the primes $q\notin S$ at which $\rho_l $ is unramified (\cite{DDT} {\bf Proposition 2.6 (3)}), the assumption implies that $\rho_{f,l} =\rho_{g,l}$ and in particular $N_f=N_g$. Now we deduce that $f=g$ from {\bf Fact 4.1}. 

\begin{flushright}
$\Box$
\end{flushright}

\begin{remark} Here is an another way to see that any $\beta$ in (18) is different from each of $\beta^{\prime}$ in (19). If necessary changing $M$ and $M^{\prime}$, let $r$ be a prime divisor of $M^{\prime}$ not dividing $M$. By the construction $\beta^{\prime}(U_r)\in {\mathbb R}$ and $\beta(U_r)\in {\mathbb C}\setminus{\mathbb R}$ and therefore $\beta$ and $\beta^{\prime}$ are different. 
\end{remark}

For a character $\alpha$ of ${\mathbb T}_0(N)$, let $S_2(\Gamma_0(N))(\alpha)$ denote the isotypic component of $\alpha$ ;
\[S_2(\Gamma_0(N))(\alpha)=\{f\in S_2(\Gamma_0(N))\,|\, f|T=\alpha(T)f,\quad \forall T\in {\mathbb T}_0(N)\}.\]

\begin{thm} (Strong multiplicity one) Let $N$ be a square free positive integer. Then there is an isomorphism as ${\mathbb T}_0(N)$-modules
\[S_2(\Gamma_0(N))=\oplus_{\alpha}S_2(\Gamma_0(N))(\alpha)\]
such that every irreducible component has dimension one and is spanned by the normalized Hecke eigenform $f_{\alpha}$. The index $\alpha$ in the decomposition runs through the set of closed points ${\rm Spec}({\mathbb T}_0(N))({\mathbb C})$ and there is an isomorphism
\[\Phi : {\mathbb T}_0(N)\otimes{\mathbb C} \simeq \prod_{\alpha\in {\rm Spec}({\mathbb T}_0(N))({\mathbb C})}{\mathbb C}\]
such that the composition with the projection $\pi_{\alpha}$ to the $\alpha$-factor is $\alpha$ :
\[\pi_{\alpha} \circ \Phi=\alpha.\]
\end{thm}
{\bf Proof.}
The previous argument and (11) show that  $S_2(\Gamma_0(N))(\alpha)$ is a ${\mathbb C}$-linear space generated by a normalized Hecke eigenform $f_{\alpha}$ and we have an irreducible decomposition of multiplicity one
\begin{equation}S_2(\Gamma_0(N))=\oplus_{\alpha}S_2(\Gamma_0(N))(\alpha).\end{equation}
The linear isomorphism
\[{\rm Hom}_{\mathbb C}({\mathbb T}_0(N), {\mathbb C})\simeq S_2(\Gamma_0(N)),\quad \rho \mapsto \sum_{m=1}^{\infty}\rho(T_m)q^m\]
implies that $\{\alpha\}$ in the right hand side of (20) is the set of closed points ${\rm Spec}({\mathbb T}_0(N))({\mathbb C})$ and $\{f_{\alpha}\}_{\alpha\in {\rm Spec}({\mathbb T}_0(N))({\mathbb C})}$ is a basis of $S_2(\Gamma_0(N))$. Now the desired decomposition of ${\mathbb T}_0(N)\otimes{\mathbb C}$ is obvious.
\begin{flushright}
$\Box$
\end{flushright}
%
Let $p$ be any prime ({\em not necessary} $p\equiv1 ({\rm mod}\,12)$) and $N$ a square free positive integer prime to $p$. We define {\em the space of $p$-new forms} $S_2(\Gamma_0(pN))_{pN/N}$ to be the orthogonal complement of $\alpha_p^{\ast}(S_2(\Gamma_0(N)))+\beta_p^{\ast}(S_2(\Gamma_0(N)))$ in $S_2(\Gamma_0(pN))$ for the Petersson inner product.
Then (11) and (14) imply
\begin{equation}S_2(\Gamma_0(pN))_{pN/N}=\oplus_{M|N}\oplus_{d|(N/M)}S_2(\Gamma_0(pM))_{new}^{(d)}\end{equation}
and by {\bf Theorem 4.1} we have a decomposition of ${\mathbb T}_0(N)$-module of multiplicity one
\begin{equation}S_2(\Gamma_0(pN))_{pN/N}=\oplus_{\chi}{\mathbb C}f_{\chi}.\end{equation}
Here $f_{\chi}$ is a normalized Hecke eigenform whose character is $\chi$. Let ${\mathbb T}_0(pN)_{pN/N}$ be the restriction of ${\mathbb T}_0(N)$ to this space. Then the set of characters in (22) coincides with ${\rm Spec}({\mathbb T}_0(pN)_{pN/N})({\mathbb C})$ and there is an isomorphism
\begin{equation} \Phi : {\mathbb T}_0(pN)_{pN/N}\otimes{\mathbb C} \simeq \prod_{\chi\in {\rm Spec}({\mathbb T}_0(pN)_{pN/N})({\mathbb C})}{\mathbb C}\end{equation}
such that the composition with the projection $\pi_{\chi}$ to $\chi$-factor is $\chi$ :
\[\pi_{\chi}\circ \Phi=\chi.\]
%
Using \cite{Ribet1990} we will clarify a relationship between $S_2(\Gamma_0(pN))_{pN/N}$ and the Ramanujan graph $G^{(l)}_{p}(N)$. \\

By the functoriality $\alpha_p$ and $\beta_p$ induce a homomorphism
\begin{equation}\alpha_p^{\ast}, \beta_p^{\ast} : J_0(N) \to J_0(pN)\end{equation}
and we define a subvariety
\[J_0(pN)_{p-old}=\alpha_p^{\ast}J_0(M)+\beta_p^{\ast}J_0(N)\subset J_0(pN)\]
which is called as {\em $p$-old subvariety}. We define {\em $p$-new subvariety} to be the quotient
\[J_0(pN)_{pN/N}=J_0(pN)/J_0(pN)_{p-old}.\]
Now we consider the actions of Hecke operators. 
Let $T$ be $T_r\, (r\nmid pN)$ or $U_l\, (l\mid N)$. Then $T$ commutes with $\alpha_p$ and $\beta_p$ and 
\begin{equation}
\begin{CD}
 J_0(N)\times J_0(N)@>\text{$\alpha_p^{\ast}\times \beta_p^{\ast}$}>>  J_0(pN)\\
 @V\text{$(T,\,T)$}VV  @V\text{$T$}VV\\
   J_0(N)\times J_0(N)  @>\text{$\alpha_p^{\ast}\times \beta_p^{\ast}$}>>   J_0(pN).
 \end{CD}
\end{equation}
and $J_0(pN)_{p-old}$ is ${\mathbb T}_0(pN)^{(p)}$-stable. By \cite{Ribet1990} {\bf Rerark 3.9}  $J_0(pN)_{p-old}$ is also preserved by $U_p$ and it is ${\mathbb T}_0(pN)={\mathbb T}_0(pN)^{(p)}[U_p]$-stable. Therefore $J_0(pN)_{pN/N}$ admits the action of ${\mathbb T}_0(pN)$ and the image of ${\mathbb T}_0(pN)$ in ${\rm End}(J_0(pN)_{pN/N})$ is temporary denoted by ${\mathbb T}^{\prime}$. Having identified the  holomorphic cotangent space of $J_0(pN)_{pN/N}$ at the origin with $S_2(\Gamma_0(pN))_{pN/N}$ let us consider the representation of ${\rm End}(J_0(pN)_{pN/N})$ on $S_2(\Gamma_0(pN))_{pN/N}$. Then the image of ${\mathbb T}^{\prime}$ in ${\rm End}(S_2(\Gamma_0(pN))_{pN/N})$ is ${\mathbb T}_0(pN)_{pN/N}$. Since
representation of ${\rm End}(J_0(pN)_{pN/N})$ on $S_2(\Gamma_0(pN))_{pN/N}$ faithful, ${\mathbb T}^{\prime}$ and ${\mathbb T}_0(pN)_{pN/N}$ are isomorphic and we identify them.\\

%
%

%
It is known that the N\'{e}ron model of $J_0(pN)_{pN/N}$ over ${\rm Spec}{\mathbb Z}$ has purely toric reduction ${\mathcal T}$ at $p$. Let us describe its character group. 
$X_0(pN)_{{\mathbb F}_p}$ has two irreducible components $Z_F$ and $Z_V$, which are isomorphic to $X_0(N)_{{\mathbb F}_p}$. Over $Z_F$ (resp. $Z_V$) the parametrized cyclic group $C_p$ of order $p$ is the kernel of the Frobenius $F$ (resp. the Verschiebung $V$).  $Z_F$ and $Z_V$ transversally intersect at enhanced supersingular points of level $N$, that is $\Sigma_N=\{{\mathbf E}_1,\cdots,{\mathbf E}_{\nu(N)}\}$. Set
\[X_N=\oplus_{i=1}^{\nu(N)}{\mathbb Z}{\mathbf E}_i\]
and we adopt $\{{\mathbf E}_1,\cdots,{\mathbf E}_{\nu(N)}\}$ as a base. We define the action of Hecke operators on $X_N$ by (7) and (8) and let ${\mathbb T}$ denote a commutative subring of ${\rm End}_{\mathbb Z}(X_N)$ generated by Hecke operators. Let us consider the boundary map of the dual graph of $X_0(pN)_{{\mathbb F}_p}$,
\[\partial : X_N \to {\mathbb Z}Z_F\oplus {\mathbb Z}Z_F, \quad \partial({\mathbf E}_i)=Z_F-Z_V.\]
Being $X^{(0)}_N$ the kernel of $\partial$, we have the exact sequence of Hecke modules
\begin{equation}
0 \to X^{(0)}_N \to X_N \stackrel{\partial}\to {\mathbb Z}\epsilon \to 0,\quad \epsilon=Z_F-Z_V.
\end{equation}
 For brevity let us write $E_i$ by $[i]$. Then 
 \[\partial([i])=\epsilon,\quad 1\leq \forall i \leq n\]
 and
\[X^{(0)}_N=\{\sum_{i=1}^{n}a_i [i]\,|\,a_i\in{\mathbb Z},\,\sum_{i=1}^{n}a_i=0\}.\]
The the restriction ${\mathbb T}_0$ of ${\mathbb T}$ to $X^{(0)}_N$ has the following description. By \cite{Ribet1990} {\bf Proposition 3.1}, $X^{(0)}_N$ is the character group of the connected component of the torus ${\mathcal T}$. By the N\'{e}ron property, ${\mathcal T}$ admits the action of ${\mathbb T}_0(pN)_{pN/N}(={\mathbb T}^{\prime})$ and the induced action on $X^{(0)}_N$ is ${\mathbb T}_0$. Therefore ${\mathbb T}_0$ is the image of ${\mathbb T}_0(pN)_{pN/N}$ in ${\rm End}_{\mathbb Z}(X^{(0)}_N)$. Since the action of ${\mathbb T}_0(pN)_{pN/N}$ on $X_N^{(0)}$ is faithful (\cite{Ribet1990} {\bf Theorem 3.10}), ${\mathbb T}_0$ and ${\mathbb T}_0(pN)_{pN/N}$ are isomorphic and they will be identified from now on.
\begin{thm} Let $N$ be a square free positive integer. There is an isomorphism as ${\mathbb T}_0(pN)_{pN/N}$-modules
\[X^{(0)}_N\otimes{\mathbb C}\simeq S_2(\Gamma_0(pN))_{pN/N}.\]
\end{thm}
{\bf Proof.}
As we have mentioned before, the action of ${\mathbb T}_0(pN)_{pN/N}$ on $X_N^{(0)}$ is faithful (\cite{Ribet1990} {\bf Theorem 3.10}). Since the characters $\{\chi\}$ in (22) are mutually different and by (23) we see every irreducible component of (22) should appear as irreducible factor of $X^{(0)}_N\otimes{\mathbb C}$. Thus $S_2(\Gamma_0(pN))_{pN/N}$ is contained in $X^{(0)}_N\otimes{\mathbb C}$. On the other hand the rank of $X^{(0)}_N$ is equal to ${\rm dim}{\mathcal T}={\rm dim}J_0(pN)_{pN/N}$. Since the holomorphic cotangent space of $J_0(pN)_{pN/N}$ at the origin is $S_2(\Gamma_0(pN))_{pN/N}$,
\[{\rm dim}X^{(0)}_N\otimes{\mathbb C}={\rm dim}S_2(\Gamma_0(pN))_{pN/N},\]
and the claim is proved. 

\begin{flushright}
$\Box$
\end{flushright}
Let us state a real version of {\bf Theorem 4.2}. Since the character of a normalized Hecke-eigen newform is real,  using (15) and (20), {\bf Theorem 4.2} yields an decomposition as a ${\mathbb T}_0(pN)_{pN/N}\otimes{\mathbb R}$-module
\[X^{(0)}_N\otimes{\mathbb R}=\oplus_{\gamma}V(\gamma),\]
where
\[V(\gamma)=\{v \in X^{(0)}_N\otimes{\mathbb R}\,|\, T(v)=\gamma(T)v\quad \forall T\in {\mathbb T}_0(pN)_{pN/N}^{(pN)}\}.\]
Here $\gamma$ is the real character of ${\mathbb T}_0(pN)_{pN/N}^{(pN)}$ which is the restriction of the character of the normalized Hecke eigen newform $f_{\gamma}$ whose level $N_{\gamma}$ satisfies
\[N_{\gamma}=pM,\quad M|N.\]
{\bf Lemma 4.1} shows that $\{\gamma\}$ are mutually different. Being $N/M=l_1\cdots l_m$ the prime decomposition, we write
\begin{equation}{\mathbb T}_0(pN)_{pN/N}\otimes{\mathbb R}=({\mathbb T}_0(pN)_{pN/N}^{(N/M)}\otimes{\mathbb R})\otimes_{\mathbb R}{\mathbb R}[U_{l_1},\cdots,U_{l_m}]\end{equation}
and $V(\gamma)$ is a ${\mathbb R}[U_{l_1},\cdots,U_{l_m}]$-module. As we have seen in the proof of {\bf Proposition 4.1}, the characteristic polynomial of $U_{l_i}$ is $P_{l_i}(U_{l_i})=U_{l_i}^2-a_{l_i}(f_{\gamma})U_{l_i}+l_i$ and ${\rm dim}_{\mathbb R}V(\gamma)=2^m$. Therefore
\[V(\gamma)\simeq {\mathbb R}[U_{l_1},\cdots,U_{l_m}]/I,\]
where  $I$ is an ideal of ${\mathbb R}[U_{l_1},\cdots,U_{l_m}]$ generated by the polynomials $\{P_{l_i}(U_{l_i})\}_{i=1,\cdots,m}$. Viewing ${\mathbb R}f_{\gamma}$ as a ${\mathbb T}_0(pN)_{pN/N}^{(N/M)}\otimes{\mathbb R}$-module, we write it by ${\mathbb R}f^{(N/M)}_{\gamma}$. Using (27) we see
\[V(\gamma)\simeq {\mathbb R}f^{(N/M)}_{\gamma} \otimes_{\mathbb R}({\mathbb R}[U_{l_1},\cdots,U_{l_m}]/I).\]
as ${\mathbb T}_0(pN)_{pN/N}\otimes{\mathbb R}$-modules.
Thus we have proved a real version of {\bf Theorem 4.2}.
\begin{thm} (Weak multiplicity one) There is an irreducible decomposition \[X^{(0)}_N\otimes{\mathbb R}=\oplus_{\gamma}V(\gamma)\]
as a ${\mathbb T}_0(pN)_{pN/N}\otimes{\mathbb R}$-module. Here $\{\gamma\}$ runs through the real characters of normalized Hecke eigen newforms $\{f_{\gamma}\}_{\gamma}$ such that the level $N_{f_{\gamma}}$ of $f_{\gamma}$ satisfies $N_{f_{\gamma}}=pM$ where $M$ is a divisor of $N$. Being $N/M=l_1\cdots l_m$ the prime decomposition, a ${\mathbb T}_0(pN)_{pN/N}\otimes{\mathbb R}$-module $V(\gamma)$ is defined to be
\[V(\gamma)\simeq {\mathbb R}f^{(N/M)}_{\gamma} \otimes_{\mathbb R}({\mathbb R}[U_{l_1},\cdots,U_{l_m}]/I).\]
Here the action of ${\mathbb T}_0(pN)_{pN/N}\otimes{\mathbb R}$ is defined via (27) and $I$ is an ideal generated by polynomials $\{P_{l_i}(U_{l_i})\}_{i=1,\cdots,m}$ where
\[P_{l_i}(U_{l_i})=U_{l_i}^2-a_{l_i}(f_{\gamma})U_{l_i}+l_i.\] 
Moreover the characters $\{\gamma\}$ are mutually different.

\end{thm}
Let $l$ be an odd prime different from $p$. Remember that $N\in {\mathcal N}_{p,l}$ is the set of square free positive integers prime to $lp$.
\begin{thm} (Monotonicity) For $N\in {\mathcal N}_{p,l}$ let $\rho^1_l(N)$ be the largest eigenvalue of the Hecke operator $T_l$ of $X^{(0)}_N\otimes{\mathbb R}$. Then for $M,N \in {\mathcal N}_{p,l}$ such that $M|N$, 
\[\rho^1_l(N) \geq \rho^1_l(M)\]
\end{thm}
{\bf Proof.} {\bf Theorem 4.2} (or {\bf Theorem 4.3}) shows that, under the decomposition (22), $\rho^1_l(N)$ is the maximum of $l$-th coefficients of Hecke eigenform $\{f_{\chi}\}_{\chi}$. By (21) we find $S_2(\Gamma_0(pM))_{pM/M}$  is contained in $S_2(\Gamma_0(pN))_{pN/N}$ and the claim is obtained.

\begin{flushright}
$\Box$
\end{flushright}

\section{Properties of the graphs}
Let $p$ be a prime satisfying $p\equiv 1 ({\rm mod}\,12)$ and $l$ be an odd prime different from $p$. Let us take $N\in {\mathcal N}_{p,l}$. For brevity we write ${\mathbf E}_i=(E_i, C_N)$ and let $\Gamma_l$ be the set of cyclic subgroups of $E_i$ of order $l$. The bijective correspondence
\[{\rm Hom}({\mathbf E}_i,\,{\mathbf E}_j)(l)/{\pm 1} \simeq \Gamma_l, \quad f \mapsto {\rm Ker}f.\]
shows that the Brandt matrix $B_p^{(l)}(N)$ is the representation matrix of $T_l$. Since $B_p^{(l)}(N)$ is symmetric, the eigenvalues are all real. It is easy to check that $\epsilon=Z_F-Z_V$ (cf. (26)) satisfies
\[T_l(\epsilon)=(l+1)\epsilon\]
and since $\partial$ in (26) commutes with $T_l$, $l+1$ is an eigenvalue of $B_p^{(l)}(N)$. Let $\delta$ be a corresponding eigenvector. Using the Eichler-Shimura relation and the Weil conjecture, {\bf Theorem 4.2} (or {\bf Theorem 4.3}) implies that the modulus of other eigenvalues are less than or equal to $2\sqrt{l}$ and 
\[X_{N}\otimes{\mathbb R}=(X^{(0)}_{N}\otimes{\mathbb R})\hat{\oplus} {\mathbb R}\delta,\]
where $\hat{\oplus}$ denotes an orthogonal direct sum.  Moreover if $N$ is generic, {\bf Theorem 4.3} and this decomposition yield a spectral decomposition  of $X_{N}\otimes{\mathbb R}$ in terms of eigenspaces of $T_l$. {\bf Theorem 4.2} implies that
\begin{equation}{\rm det}[1-B_p^{(l)}(N)t+lt^2]=(1-t)(1-lt){\rm det}[1-T_lt+lt^2 | S_2(\Gamma_0(pN))_{pN/N}].\end{equation}

\begin{thm} For any $N\in {\mathcal N}_{p,l}$,  $G_p^{(l)}(N)$ is a connected regular Ramanujan graph of degree $l+1$ not bipartite.
\end{thm}
{\bf Proof.} By construction $G_p^{(l)}(N)$ is a regular graph of degree $l+1$. Let us investigate the eigenvalues of the adjacency matrix $B_p^{(l)}(N)$. As we have seen, $l+1$ is an eigenvalue of $B_p^{(l)}(N)$ and the modulus of other eigenvalues are less than or equal to $2\sqrt{l}$. Thus $G_p^{(l)}(N)$ is a Ramanujan graph.
By the equation (1) (see also (2)), $0$ is an eigenvalue of the Laplacian with multiplicity one and we see that $G_p^{(l)}(N)$ is connected. In general a connected finite regular graph of degree $d$ is bipartite if and only if $\pm d$ are eigenvalues of the adjacency matrix (\cite{Terras}). Therefore $G_p^{(l)}(N)$ is not bipartite.
\begin{flushright}
$\Box$
\end{flushright}
Now {\bf Theorem 1.1} is a direct consequence of the equation (1) (see also (2)), {\bf Theorem 4.4} and {\bf Theorem 5.1}.\\

\noindent
{\bf Proof of Theorem 1.2} 
Set $N=q$ and we use the decomposition (21). Since $S_2(\Gamma_0(p))^{(q)}$ is isomorphic to $S_2(\Gamma_0(p))$ as a ${\mathbb T}_0(pq)^{(pq)}$-module, we see
\[S_2(\Gamma_0(pq))_{pq/q}=S_2(\Gamma_0(pq))_{new}\oplus S_2(\Gamma_0(p))^{\oplus 2}\]
as ${\mathbb T}_0(pq)^{(pq)}$-modules and 
%
\[\frac{{\rm det}(1-B_p^{(l)}(q)t+lt^2)}{{\rm det}(1-B_p^{(l)}(1)t+lt^2)^2}=\frac{{\rm det}(1-T_lt+lt^2\mid S_2(\Gamma_0(pq))_{new})}{(1-t)(1-lt)}=\frac{{\rm det}(1-B_q^{(l)}(p)t+lt^2)}{{\rm det}(1-B_q^{(l)}(1)t+lt^2)^2}\]
by (28). On the other hand {\bf Fact 2.2} implies,
\[ \chi(G^{(l)}_p(q))-2\chi(G^{(l)}_p(1))=\frac{(p-1)(q-1)(1-l)}{24}=\chi(G^{(l)}_q(p))-2\chi(G^{(l)}_q(1))\]
and the claim follows from {\bf Fact 2.1}. 

\begin{flushright}
$\Box$
\end{flushright}

\noindent
{\bf Proof of Theorem 1.3}
Let us recall the decomposition (22)
\[S_2(\Gamma_0(pN))_{pN/N}=\oplus_{\chi}{\mathbb C}f_{\chi},\]
where $f_{\chi}$ is a normalized Hecke eigenform. Then the second largest eigenvalue $\rho^1_l(N)$ of $B_p^{(l)}(N)$ is the maximum of $\{a_l(f_{\chi})\}_{\chi}$ by {\bf Theorem 4.2} and satisfies $\rho^1_l(N)\leq 2\sqrt{l}$ by {\bf Theorem 5.1}. Let $\{r_i\}_{i=1}^{\infty}$ be the set of primes and  $N_k=\prod_{i=1}^{k}r_i$. Then by {\bf Theorem 4.4},  $\rho^1_l(N_k)$ is monotone increasing   for $k$. In general let $\{G_i\}_{i}$ be an infinite family of connected $d$-regular graphs satisfying
\[\lim_{i\to \infty}|V(G_{i})|=\infty.\]
Then it is known that
\[\liminf_{i\to \infty}\rho^1(G_i)\geq 2\sqrt{d-1}\]
by Alon and Boppana (\cite{Alon} \cite{Alon-Milman}). We will use this fact. Since $\{G_p^{(l)}(N_k)\}_k$ is an infinite family of connected regular Ramanujan graphs of degree $l+1$ with
\[\lim_{k\to\infty}|V(G_p^{(l)}(N_k))|=\lim_{k\to \infty}\frac{(p-1)\prod_{i=1}^{k}(1+r_i)}{12}=\infty,\]
we see
\[\lim_{k\to \infty} \rho^1_l(N_k)=2\sqrt{l},\]
and 
\[\lim_{k\to \infty}{\rm Max}\{a_l(f_{\chi})\,:\, S_2(\Gamma_0(pN_k))_{pN_k/N_k}=\oplus_{\chi}{\mathbb C}f_{\chi}\}=2\sqrt{l}.\]
Since $S_2(\Gamma_0(pN_k))_{pN_k/N_k}$ is a subspace of $S_2(\Gamma_0(pN_k))$,  the remaining claim immediately follows from this result and the decomposition in {\bf Theorem 4.1}.
\begin{flushright}
$\Box$
\end{flushright}
The proof implies the following corollary.
\begin{cor} Let $p$ be a prime satisfying $p\equiv 1 ({\rm mod}\,12)$ and $l$ an odd prime with $l\neq p$. Then for any set of mutually distinct primes $\{r_i\}_{i=1}^{\infty}$ which are different from $l$ and $p$, there is a sequence of normalized Hecke eigenforms $\{f_i\}_i$ of weight $2$ such that $f_i \in S_2(\Gamma_0(pr_1\cdots r_i))_{new}$ and 
\[\lim_{i\to\infty}a_l(f_i)=2\sqrt{l}.\]

\end{cor}

%

\noindent

\end{document}